\documentclass[11pt]{amsart}
\usepackage{amsmath}
\usepackage{amsthm}
\usepackage{amsopn}
\usepackage{amssymb}
\usepackage{fullpage}
\usepackage{enumerate}
\numberwithin{equation}{section}
\usepackage{float}
\usepackage{graphicx}
\usepackage{titlesec}
\usepackage{dsfont}
\usepackage{fancyhdr}
\usepackage{textcomp}
\usepackage{fancyhdr}
\newcommand\smallO{
  \mathchoice
    {{\scriptstyle\mathcal{O}}}
    {{\scriptstyle\mathcal{O}}}
    {{\scriptscriptstyle\mathcal{O}}}
    {\scalebox{.7}{$\scriptscriptstyle\mathcal{O}$}}
  }

\newtheorem{lemma}{Lemma}[section]
\newtheorem{theorem}{Theorem}[section]
\newtheorem{proposition}{Proposition}[section]

\theoremstyle{definition}
\newtheorem{remark}{Remark}[section]

\DeclareMathOperator{\Capac}{cap}

\DeclareMathOperator{\dist}{dist}
\DeclareMathOperator{\arctanh}{arctanh}

\DeclareMathOperator{\Ar}{A}
\DeclareMathOperator{\bs}{ \Omega^{\star}}
\DeclareMathOperator{\Le}{L}

\DeclareMathOperator{\T}{\mathbb{T}}
\DeclareMathOperator{\supp}{supp}
\DeclareMathOperator{\caph}{caph}

\DeclareMathOperator{\D}{\mathbb{D}}
\DeclareMathOperator{\RE}{Re}
\DeclareMathOperator{\IM}{Im}
\titleformat{\section}
  {\normalfont\scshape}{\thesection}{1em}{\centering}
\titleformat{\subsection}[runin]
  {\bfseries}{\thesubsection}{1em}{}

\begin{document}
\title{On the spectral value of Semigroups of Holomorphic Functions}


\author{Maria Kourou}  

\address{Department of Mathematics, Aristotle University of Thessaloniki, 54124, Thessaloniki, Greece}
\email{mkouroue@math.auth.gr}

\fancyhf{}
\renewcommand{\headrulewidth}{0pt}
\fancyhead[RO,LE]{\small \thepage}
\fancyhead[CE]{\footnotesize MARIA KOUROU}
\fancyhead[CO]{\footnotesize On the Spectral Value of Semigroups of Holomorphic Functions} 

\fancyfoot[L,R,C]{}
\subjclass[2010]{Primary 31A15, 30D05, 47D06; Secondary 30C20, 30C85, 37C10}

\date{}
\keywords{Semigroup of holomorphic functions, harmonic measure, extremal distance, condenser capacity}

\begin{abstract}
Let $(\phi_t)_{t \geq 0}$ be a semigroup of holomorphic self-maps of the unit disk $\D$ with Denjoy-Wolff point $\tau=1$. The angular derivative is $\phi_t^{\prime}(1)= e^{-\lambda t}$, where $\lambda \geq 0$ is the spectral value of $(\phi_t)$. If $\lambda>0$ the semigroup is hyperbolic, otherwise it is parabolic.
Suppose $K$ is a compact non-polar subset of $\D$ with positive logarithmic capacity. We specify the type of the semigroup by examining the asymptotic behavior of $\phi_t(K)$. We provide a representation of the spectral value of the semigroup with the use of several potential theoretic quantities e.g. harmonic measure, Green function, extremal length, condenser capacity. 



%
\end{abstract}

\maketitle

\section{Introduction}

A one-parameter semigroup is a family $(\phi_t)_{t \geq 0}$ of holomorphic functions in $\D$, where 
\begin{enumerate}[(i)]
\item $\phi_0$ is the identity map;
\item $\phi_{t+s}(z) = \phi_t \left( \phi_s (z) \right)$, for every $t,s \geq 0$ and $z \in \D$;
\item $\phi_t(z) \xrightarrow{t \to 0^{+}} z$, uniformly on compacta in $\D$.
\end{enumerate}

The introduction of the theory of one-parameter semigroups was made by Berkson and Porta in \cite{berksonporta}, where their strong connection with composition operators is emphasized. Later, one-parameter semigroups have been an issue of scientific interest, especially due to their application in other mathematical areas  such as dynamical systems, Markov stochastic processes, etc. 
Recent advances in this area, as well as, the classical theory of semigroups can be found in \cite{BetsDescript, Booksem, AnalyticFlows, shoikhetelin, kelgiannis, shoikhet}. 

By the continuous form of Denjoy-Wolff theorem, for every such one-parameter semigroup, there exists a unique attractive fixed point $\tau \in \overline{\D} $, which is called the \textit{Denjoy-Wolff point} of $(\phi_t)$; see \cite[Theorem 1.4.17]{abate}. 
Fix $z \in \D$. Let's denote by $ \gamma_z$ the curve 
$ \gamma_z: [0, + \infty)  \to \D$ with $\gamma_z(t)= \phi_t(z).$
This curve is called the \textit{trajectory} of $z$ and 
\begin{equation}\label{limtraj}
\lim_{t \to + \infty} \gamma_z(t) =\lim_{t \to +\infty} \phi_t(z)= \tau.
\end{equation}

Depending on the position of $\tau$, we have the following classification of one-parameter semigroups.  
Supposing that $\phi_t$ is not an elliptic automorphism of $\D$ for any $t \geq 0$, the semigroup $(\phi_t)$ is called \textit{elliptic}, if $\tau$ lies inside the unit disk. 

In the present work, we are interested in the case where $(\phi_t)$ is a non-elliptic semigroup. Hence the Denjoy-Wolff point lies on the unit circle $\T$ and without loss of generality, we can assume that it is equal to $1$. 

The angular derivative $\phi_t^{\prime}$ at $1$ is given by the limit
$$\phi_t^{\prime}(1):= \angle \lim_{z \to 1} \frac{\phi_t(z)-1 }{z-1} =e^{- \lambda t}\leq 1,$$
where the number $\lambda$ is the \textit{spectral value} of the semigroup. 
If $\lambda>0$, the semigroup $(\phi_t)$ is called \textit{hyperbolic}, whereas, if $\lambda =0$, the semigroup $(\phi_t)$ is called \textit{parabolic}. 

In \cite[\S\ 9.1, Theorem 9.1.2]{Booksem}, it was proved that
\begin{equation}\label{divrate}
\frac{\lambda}{2} = \lim_{t \to + \infty} \frac{d_{\D}(z, \phi_t(z))}{t}, \quad z \in \D
\end{equation} 
where $d_{\D}$ denotes the hyperbolic distance in the unit disk.

The aim of the present work is to provide a characterization of the spectral value using the trajectory of a compact set instead of the trajectory of a point under the semigroup. 
In order to achieve this, we use a variety of potential theoretic quantities such as the harmonic measure, the extremal length and the condenser capacity. 

Suppose $K$ is a compact non-polar subset of $\D$. This means that its logarithmic capacity $\Capac K$ is strictly positive; see Subsection \ref{diameter}. Then $\phi_t(K)$ approaches the unit circle, as $t$ increases, and asymptotically, it shrinks to the Denjoy-Wolff point of the semigroup, due to \eqref{limtraj}. 

For a point $z \in \D$, the harmonic measure $\omega (z, \phi_t (K), \D)$ is the Perron solution of the Dirichlet problem on $\D \setminus \phi_t(K)$ with boundary values $1$ on the boundary of $\phi_t(K)$ and $0$ on $\T$, as one can see in detail in Subsection \ref{harmonicmeasure}.
Using Beurling's estimate \cite[Theorem H.7]{margarnett}, as the reader may observe in detail in Section \ref{lemmas}, for every $z \in K$, there exists a constant $C$ such that 
$$ \omega(z, \phi_t(K), \D) \leq C e^{- \pi \lambda_{\D}(K, \phi_t(K)) },$$
where $\lambda_{\D}(K, \phi_t(K))$ denotes the extremal distance between the compact sets $K$ and $\phi_t(K)$.
The extremal distance is the most natural way to measure the distance between sets inside a domain in a conformally invariant way. 
We can say that it is in a way equivalent to the hyperbolic distance, when it comes to compact sets. For the definition and the properties of extremal distance, we refer to Subsection \ref{extremallength}.  

The reciprocal of the extremal distance between $K$ and $\phi_t(K)$ can be expressed as the capacity of the generalized condenser $(\D, K, \phi_t(K))$ on the premise that both $K$ and $\phi_t(K)$ are non-polar; see Subsection \ref{extremallength}. 
The capacity of $(\D, K, \phi_t(K))$ is defined as 
$$ \Capac(\D, K, \phi_t(K)) = \int_{\D \setminus (K \cup \phi_t(K))} | \nabla u|^2 dA,$$
where $A$ is the Lebesgue measure on the complex plane and $u$ has boundary values $0$ on $\partial K$, $1$ on $\partial \phi_t(K)$ and $\frac{\partial u}{\partial n}=0$ on $\T$. 
More information on generalized condensers follows in Section \ref{condensers}. 
  
The question that we deal with is the inverse problem of the one addressed in the articles \cite{bkkp2} and \cite{kourou3}. Knowing how $\phi_t(K)$ behaves, as $t \to +\infty$, can we distinguish the type of the semigroup? Moreover, can we obtain an analogous form of the spectral value of a semigroup using compact sets instead of points? 
In fact, we reach the following outcome. 

\begin{theorem}\label{main}
Let $(\phi_t)$ be a non-elliptic semigroup of holomorphic self-maps of $\D$ and let $K$ be a compact non-polar subset of $\D$. Then, as $ t \to + \infty$, 
\begin{enumerate}[\normalfont(i)]
\item $2 t^{-1} d_{\D}(z,\phi_t(w))\rightarrow \lambda$, for every $z,w \in \D$
\item $-t^{-1} \log g_{\D}(z, \phi_t(w)) \rightarrow   \lambda$, for every $z, w \in \D$ 
\item $-t^{-1} \log \omega (z, \phi_t(K), \D) \rightarrow  \lambda$, for every $z \in K$
\item $ \pi t^{-1} \lambda_{\D} (K, \phi_t(K))  \rightarrow  \lambda$
\item $\frac{1}{\pi} t \Capac(\D, K, \phi_t(K)) \rightarrow \frac{1}{\lambda}$.
\end{enumerate}
\end{theorem}

The spectral value, as seen in \eqref{divrate}, depends on the hyperbolic distance between a point and its image under $\phi_t$. We prove in Theorem \ref{main} that in the case where arbitrary compact non-polar sets are examined, the hyperbolic distance is substituted by the extremal distance in the unit disk. 

To continue with, the class of parabolic semigroups is further partitioned into two categories with the use of the \textit{hyperbolic step}; \cite[Definition 9.1.5]{Booksem}. Let $u>0$.  
The hyperbolic step of order $u$ of $(\phi_t)$ at $z \in \D$ is defined as 
\begin{equation}\label{hyperbstep}
s_u(\phi_t,z):= \lim_{r \to + \infty} d_{\D}( \phi_r(z), \phi_{r+u}(z)). 
\end{equation}

If $s_u(\phi_t,z)=0$, for all $u>0$ and $z \in \D$, the parabolic semigroup is of \textit{zero hyperbolic step}. 
If $s_u(\phi_t,z)>0$, for some $u>0$ or $z \in \D$, then the parabolic semigroup is of \textit{positive hyperbolic step}.
There is an immediate connection between the spectral value and the hyperbolic step:
\begin{equation}\label{connection}
 \frac{\lambda}{2} = \lim_{u \to + \infty} \frac{s_u(\phi_t,z)}{u},
 \end{equation} 
as we can see in \cite[Proposition 9.1.7]{Booksem}. 

Combining Theorem \ref{main} with \eqref{connection}, we obtain that the extremal distance $\lambda_{\D}(K, \phi_t(K))$ and the hyperbolic step $s_u(\phi_t,z)$ have similar asymptotic behavior.

\begin{proposition}\label{asymptoticintegral}
As $u \to +\infty$, $$s_u( \phi_t, z) - \frac{\pi}{2} \lambda_{\D}(K, \phi_u(K)) = \smallO(u), \quad z \in K. $$ 
\end{proposition}

The article is structured in the following way. In Section \ref{background}, we provide some information concerning the tools that are used in the proofs. They mainly consist of conformally invariant quantities from potential theory. Afterwards, in Section \ref{lemmas}, some auxiliary lemmas are presented and proved, whereas, Theorem \ref{main} is proved in Section \ref{proof}. In addition, Proposition \ref{asymptoticintegral} is proved in Section \ref{proofprop}, where the reader may also find a connection between the extremal distance and the hyperbolic step of sufficiently large order.

\section{Background Material}\label{background}

\subsection{Koenigs function of a Semigroup}

For every semigroup $(\phi_t)$ with Denjoy-Wolff point $1$, there exists a conformal function $h$ that maps $\D$ onto a simply connected domain $\Omega$, such that 
\begin{equation}\label{mapofh}
h(\phi_t(z) )= h(z)+t, \quad \forall z \in \D, \, \forall t \geq 0.
\end{equation}
This mapping is called \textit{Koenigs function} and its main property is to linearize the trajectories of the points in $\D$.
The simply connected domain $\Omega= h(\D)$ is \textit{convex in the horizontal direction}, as $\{ w+s: s>0 \} \subset \Omega$, for every $w \in \Omega$. This occurs since the trajectory of each point $z \in \D$ is mapped onto a half-line; see \eqref{mapofh}. 

From \cite{AnalyticFlows}, the semigroup is hyperbolic if and only if $\Omega$ is contained in a horizontal strip.
Furthermore, according to \cite[Theorem 1]{BetsDescript}, if $(\phi_t)$ is parabolic of positive hyperbolic step, then $\Omega $ is contained in a horizontal half-plane. In both cases, we can find the smallest horizontal domain that contains $\Omega$, which is called the \textit{base space} of the semigroup.
The base space of a parabolic semigroup of zero hyperbolic step is the whole complex plane. 
We denote by $\bs$ the associated base space of the semigroup, regardless of the type of the semigroup. 

 \subsection{Generalized Condensers}\label{condensers}

Suppose $\Omega$ is a domain in $\mathbb{C}$.  The triple $\mathcal{C} := (\Omega, \{E_i \}_{i=1,2}, \{0,1\} )$, where $E_1$ and $E_2$ are closed non-polar pairwise disjoint sets in $\Omega$, is called \textit{generalized condenser}. 
The set $\overline{\Omega} \setminus (E_1 \cup E_2)$ is called \textit{field of the condenser}, whereas $E_i$, $i=1,2$, are called \textit{plates of the condenser}.
In order to simplify the above notation, we denote the condenser by $\mathcal{C} := ( \Omega , E_1, E_2)$. The capacity of $\mathcal{C}$ is 
$$ \Capac \mathcal{C} : = \int_{\Omega \setminus (E_1 \cup E_2)} | \nabla u|^2 dA,$$
where $A$ is the Lebesgue measure on the complex plane and $u$ is the \textit{equilibrium potential} of the condenser. More specifically, $u$ is the solution of the mixed Dirichlet-Neumann problem, with boundary values $0$ on $\partial E_1$, $u=1$ on $\partial E_2$ and $\frac{\partial u}{\partial n}=0$ on $\partial \Omega$.

An important property of condenser capacity is conformal invariance. If $f$ is a conformal mapping on the field of a condenser $\mathcal{C}$, then 
\begin{equation}\label{confinvcapacity}
\Capac \mathcal{C} = \Capac f(\mathcal{C}). 
\end{equation}

More information on condenser capacity can be found in \cite{dubininbook}.

\subsection{Euclidean $n$-th Diameter - Logarithmic capacity} \label{diameter}
Let $K$ be a compact subset of $\mathbb{C}$. The euclidean $n$-th diameter of $K$ is 
\begin{equation}\label{eudiameter}
d_n (K) = \sup_{w_{\mu}, w_{\nu} \in K} \prod_{1 \leq \mu < \nu \leq n} |w_{\mu} - w_{\nu}|^{\frac{2}{n(n-1)}}
\end{equation}
and the supremum is attained, since $K$ is compact, for a $n$-tuple of points, which is called \textit{Fekete $n$-tuple} and it may not be unique; see \cite[Definition 5.5.1]{ransford}. 
The \textit{logarithmic capacity} of $K$ is the limit
$$\Capac K= \lim_{ n \to + \infty} d_n(K). $$  

If two sets $A,B \subset \mathbb{C}$ differ on a set of zero logarithmic capacity (meaning $\Capac( A \setminus B)= \Capac (B \setminus A)=0$) we say that they are nearly everywhere (n.e.) equal. Sets of zero logarithmic capacity are called \textit{polar sets} and they are negligible in the view of potential theory.

\subsection{Green function-Green capacity}

Let $D$ be a domain of the extended complex plane $\widehat{\mathbb{C}}$. 
The \textit{Green function} of $D$ with \textit{pole} at $w \in D$  is denoted by $g_D(\cdot, w)$ and satisfies the following conditions:
\begin{enumerate}[(i)]
\item $g_D(\cdot,w) $ is positive harmonic on $D \setminus \{w \}$ and bounded outside every neighborhood of $w$
\item $g_D(\cdot,w)+ \log| \cdot -w|$ is harmonic on $D$ 
\item $g_D(w,w)=\infty$, and for $z \to w$
$$ g_D(z,w) =   \begin{cases}
 \log |z| + \mathcal{O} (1), \, w=\infty \\
	-\log|z-w| + \mathcal{O}(1), \, w \neq \infty
  \end{cases} $$
\item $g_D(\cdot,w) = 0$ on the boundary $ \partial D$. 
\end{enumerate}

For instance, the Green function of the unit disk $\D$ is equal to  
\begin{equation}\label{greenformula}
g_{\D}(z,w) = \log \left| \frac{1 - z \overline{w}}{z-w} \right|,
\end{equation}
for $z, w \in \D$; see \cite[p.109]{ransford}. 
If the boundary of a domain $D$ is non-polar, then the Green function $g_{D}$ exists and it is unique. In this case, the domain $D$ is called \textit{Greenian}. The Green function is symmetric,
for every $z,w \in D$. Moreover, the Green function is conformally invariant; e.g. \cite[ Theorem 4.4.4]{ransford}.
For each measure $\mu$ with support in $D$, we can define its \textit{Green potential}
$$G_{\mu}^D(x) = \int g_D(x,y) d \mu (y) ,$$
for $x \in D$.
The Green potential is a superharmonic function on $D$ and harmonic on $D \setminus \supp \{ \mu \}$.  

Let $E$ be a subset of $D$. The \textit{Green energy} of $E$ with respect to $D$ is defined as 
$$ V(E,D) : =  \inf_{\mu} \int G_{\mu}^D (y) d \mu(y),$$
where the infimum is taken over all Borel measures $\mu$ with compact support in $E$ and $\mu(E)=1$. 
If $E$ is compact, the infimum is attained for a Borel measure $\mu$, which is called \textit{Green equilibrium measure} of $E$. 
For further information on Green potentials, we refer to \cite{margarnett, helms}. 

\subsection{Extremal distance}\label{extremallength}

Another conformally invariant quantity that has a strong relation with condenser capacity is extremal length. Let $\{ \Gamma \}$ be a family of rectifiable curves in a domain $\Omega \subset \mathbb{C}$.  
Consider all the non-negative Borel measurable functions $\rho(z)$ on $\Omega$ such that the area
$$\Ar ( \rho, \Omega) = \iint_{\Omega} \rho(z)^2 dxdy \in (0, + \infty).$$
The $\rho$-length of $\{ \Gamma\}$ is 
$$\Le \left(\rho, \{ \Gamma \} \right) = \inf_{ \gamma \in \{\Gamma \} } \int_{\gamma} \rho(z) |dz|. $$

The \textit{extremal length} of $\{ \Gamma \}$ is defined as 
$$\lambda_{\Omega} ( \{ \Gamma \}) = \sup_{\rho} \frac{\Le (\rho, \{\Gamma \})^2}{\Ar (\rho, \Omega)}, $$
where the supremum is taken over all the functions $\rho$ with $0< \Ar (\rho, \Omega) < + \infty$.
The conformal invariance of extremal length is a direct aspect of the definition, as we can observe in \cite[p.130]{margarnett}.
The most important instance of extremal length is the \textit{extremal distance} between sets in a domain. Let $E,F \subset \overline{\Omega}$. The extremal distance from $E$ to $F$ is 
\begin{equation}\label{extconj}
\lambda_{\Omega} (E,F) : = \lambda_{\Omega \setminus (E \cup F)} (E,F) = \lambda_{\Omega} (\{ \Gamma\}),
\end{equation}
where $\{\Gamma \}$ is the family of arcs or closed curves in $\Omega$ that join $E$ and $F$. 
The \textit{conjugate extremal distance} between $E$ and $F$ is defined as $\lambda_{\Omega}^{\star}(E,F):= \lambda_{\Omega} ( \{ \Gamma^{\star} \}),$
where $\{ \Gamma^{\star} \}$ is the family of all curves that separate $E$ from $F$.
As we can see in \cite[p.491]{margarnett}, the conjugate and extremal distance satisfy the following relation
$$\lambda_{\Omega}(E,F) =\frac{ 1}{ \lambda_{\Omega}^{\star}(E,F)}.$$

Two basic examples of calculation of the extremal length are the cases of a rectangle and an annulus; see \cite[p.131-3]{margarnett}. 
We will need the following two properties of extremal length. 
\begin{lemma}[The extension rule]\cite[Theorem 4-1]{ahlforsconformal}\label{extension}
Let $\Omega$ be a domain in $\mathbb{C}$. If $\{ \Gamma \} \subset \{ \Gamma^{\prime} \}$ or equivalently, every $\gamma \in \{ \Gamma \}$ contains some $\gamma^{\prime} \in \{ \Gamma^{\prime} \}$, then $$\lambda_{\Omega} (\{ \Gamma^{\prime} \}) \leq \lambda_{\Omega} (\{ \Gamma \}).$$
\end{lemma}

\begin{lemma}\cite[Theorem 2.8]{ohtsuka}\label{non-disjointparallel}
Let $\{ \Gamma_1 \},...,\{ \Gamma_n \}$ be families of rectifiable curves in a domain $\Omega \subset \mathbb{C}$. Then 
$$\lambda_{\Omega} \left( \bigcup_{k=1}^n \{\Gamma_k \} \right)^{-1} \leq \sum_{k=1}^{n} \frac{1}{\lambda_{\Omega} (\{\Gamma_k \})}. $$

\end{lemma}

\begin{lemma}\cite[Theorem 4-5]{ahlforsconformal}
Suppose $K,E$ are two compact non-polar subsets of $\Omega$. The extremal distance $\lambda_{\Omega}(K,E)$ is the reciprocal of the Dirichlet integral $$D(u)= \iint_{\Omega} | \nabla u|^2 dA,$$
where $u$ solves the following mixed Dirichlet-Neumann problem:
\begin{itemize}
\item $\Delta u=0$, on $\Omega \setminus (E \cup K)$
\item $u=0$, on $\partial K$ and $u=1$, on $\partial E$
\item $\frac{\partial u}{\partial n}=0$, on $\partial \Omega$.
\end{itemize}
\end{lemma}

\begin{remark}
According to Section \ref{condensers}, the above Dirichlet integral coincides with the capacity of the generalized condenser $(\Omega, K,E)$. Hence
\begin{equation}\label{extdistcap}
\lambda_{\Omega}(K,E)=\frac{1}{\Capac(\Omega,K,E)}.
\end{equation}
\end{remark}

For further information on extremal length the reader may refer to \cite[Chapter 4]{ahlforsconformal}, \cite[Chapter IV]{margarnett} and \cite{ohtsuka}.

\subsection{Harmonic Measure}\label{harmonicmeasure}

Let $D$ be a proper subdomain of $\widehat{\mathbb{C}}$ with non-polar boundary and $\mathcal{B}(\partial D)$ be the $\sigma$-algebra of all Borel sets of $\partial D$.   
Suppose $E \in \mathcal{B}(\partial D)$. The harmonic measure of $E$ at a point $z \in D$ is the solution of the generalized Dirichlet problem in $D$ with boundary values $1$ on $E$ and $0$ on $\partial D \setminus E$. 
Moreover, the harmonic measure $ \omega (z, E, D)$ is the supremum of all subharmonic functions $u$ on $D$ that satisfy $$\limsup_{w \to \zeta} u(w) \leq \mathds{1}_E (\zeta), \quad \forall \zeta \in \partial \D. $$

For a fixed $E \in \mathcal{B}(\partial D)$, $\omega ( \cdot , E, D)$ is a harmonic and bounded function on $D$. 
In addition, for a fixed point $z \in D$, the map 
$ \omega ( z, \cdot, D) : \mathcal{B}(\partial D)  \to [0,1]$ with $E \mapsto \omega ( z, E , D)$
is a Borel probability measure on $\partial D$. 
In addition, if $\zeta$ is a regular point of $\partial D$, which lies outside the relative boundary of $E$ in $\partial D$, then
$ \omega ( z, E, D) \xrightarrow{ z \to \zeta} \mathds{1}_E (\zeta)$. 

Furthermore, the harmonic measure, as a probability measure, has an additional significant interpretation. Suppose $D$ is a domain on $\mathbb{C}$ and $E$ a Borel subset of $\partial D$. Let $B_t$, $t>0$, be a Brownian motion on the complex plane starting from a point $z \in D$. Let $t_0= \inf \{ t>0 : B_t \notin D \setminus E \}$ be the first exit time of $B_t$ from $ D \setminus E$. The harmonic measure $\omega(z, E, D)$ is the probability of $B_{t_0} \in E$. 

A major property of the harmonic measure is conformal invariance; see \cite[\S\ 4.3]{ransford}.
For the sake of simplicity, if $E$ is a compact subset of $D$ with positive logarithmic capacity, we use the notation
$\omega(z, E, D) := \omega(z, \partial E, D \setminus E)$. 
The following inequality provides us with a connection between harmonic measure and extremal length. 
\begin{lemma}[Beurling's estimate]\cite[Theorem H.7]{margarnett}\label{h7}
Let $D$ be a domain and let $E$ be a finite union of arcs contained in one component $C_1$ of $\partial D$. Then 
$$\omega(z_0, E, D) \leq \frac{8}{\pi} \exp \left\{-\pi \sup_{\sigma} \lambda_{D} (\sigma,E) \right\}, $$
where the supremum is taken over all Jordan arcs $\sigma$ in $\mathbb{C}$ connecting $z_0 \in D$ with $C_1 \setminus E$.
\end{lemma}

\begin{remark}
Let's denote by $\{ \Gamma_1 \}$ the family of curves that join $E$ with the arcs $\sigma$, which connect $z_0 \in D$ with $C_1 \setminus E$, as in Lemma \ref{h7}. Moreover, let's denote by $\{ \Gamma_2 \}$ the family of curves that join $E$ with the Jordan arcs $\sigma^{\prime}$, which connect $z_0$ with $\partial D \setminus E$. 
Obviously, $\{ \Gamma_1 \} \subset \{ \Gamma_2 \}$ and according to Lemma \ref{extension}, $\lambda ( \{ \Gamma_1 \})\geq \lambda ( \{ \Gamma_2 \})$, which leads to $\sup_{\sigma} \lambda_{D} ( \sigma, E) \geq \sup_{\sigma^{\prime}} \lambda_D (\sigma^{\prime},E)$. 
Thus, we obtain the following inequality 
\begin{equation} \label{generh7}
\omega(z_0, E, D) \leq \frac{8}{\pi} \exp \left\{-\pi \sup_{\sigma^{\prime}} \lambda_{D} (\sigma^{\prime},E) \right\} ,
\end{equation}
where the supremum is taken over all Jordan arcs $\sigma^{\prime}$ connecting $z_0$ to $\partial D \setminus E$. 
Hence \eqref{generh7} can substitute the inequality in Lemma \ref{h7}. 
\end{remark}

\begin{proposition}[Markov Property for harmonic measure]\cite[p.88]{portstone}\label{markov}
Suppose $\Omega$ is a domain on $\mathbb{C}$. Let $S$ be a subdomain of $\Omega$, $z \in S$ and $A$ a Borel subset of $\partial \Omega$. Then 
$$\omega(z,E,\Omega) =\int_{\partial S} \omega (z, da ,S) \cdot \omega(a,E, \Omega). $$
\end{proposition} 

Information and detailed theory of harmonic measure can be found in \cite{ahlforsconformal, margarnett, ohtsuka} and \cite[Chapter 4]{ransford}. 


\subsection{ Hyperbolic metric - Hyperbolic Capacity}\label{hypmetric}

The density of the hyperbolic metric in $\D$ is $ \lambda_{\D} (z)= (1-|z|^2)^{-1}$.
Suppose $f : \D \to U$ is a conformal mapping, where $U$ is a simply connected domain of $\mathbb{C}$.
The hyperbolic density on $U$ is equal to
\begin{equation}\label{confdens}
\lambda_U (f(z)) = \frac{ \lambda_{\D} (z)}{\left| f^{\prime}(z) \right| }.
\end{equation}

\begin{lemma}[Schwarz-Pick Lemma]\cite[Theorem 6.4]{bearmin}
Suppose $\Omega_1, \Omega_2$ are simply connected proper subregions of $\mathbb{C}$ and that $f: \Omega_1 \to \Omega_2$ is holomorphic. Then 
$$\lambda_{\Omega_2}(f(z)) | f^{\prime}(z)| \leq \lambda_{\Omega_1}(z),$$
with equality holding if and only if $f$ is conformal. 
\end{lemma}

The above Lemma easily implies that for $\Omega_1 \subset \Omega_2$, it holds
\begin{equation}\label{subhypden}
\lambda_{\Omega_2}(z) \leq \lambda_{\Omega_1}(z), 
\end{equation}
for every $z \in \Omega_1$. 
The hyperbolic distance in the unit disk, for $z,w \in \D$, 
$$d_{\D}(z,w)= \inf_{\gamma} \int_{\gamma} \lambda_{\D}(z) |dz|, $$
where the infimum is taken over all curves $\gamma$ that join $z,w$ in $\D$. 
In fact, it is equal to
\begin{equation}\label{hyperbolicdistance}
d_{\D} (z,w)= \arctanh \rho_{\D}(z,w),
\end{equation}
where $$\rho_{\D}(z,w) = \left| \frac{z-w}{1- \bar{z}w} \right| $$ is the pseudo-hyperbolic distance in the unit disk $\D$. 
Hyperbolic distance is invariant under conformal mappings. Suppose $U$ is a simply connected domain of $\mathbb{C}$ and $f : \D \to U$ is conformal. Then $d_{\D} (z,w) = d_U (f(z),f(w)),$
for every choice of $z,w \in \D$.  
Furthermore, for $z, w \in \D$, it is true that
\begin{equation}\label{greendist}
g_{\D} (z, w) = -\log \tanh d_{\D}(z,w)= - \log \rho_{\D}(z,w) .
\end{equation}

Moreover, for a compact set $K \subset \D$, its hyperbolic $n$-th diameter is defined as 
\begin{equation}\label{hypdiameter}
d_{n,h}^{\D} (K)  = \sup_{w_{\mu}, w_{\nu} \in K} \prod_{1 \leq \mu < \nu \leq n} 
\rho_{\D}(w_{\mu}, w_{\nu})^{\frac{2}{n(n-1)}},
\end{equation}
where the supremum is attained for a $n$-tuple of points that are distinct and lie on the boundary. This is called \textit{Fekete $n$-tuple} and in general, it may not be unique; see \cite[Theorem 1.22]{dubininbook}. 
The \textit{hyperbolic capacity} of $K$ is 
\begin{equation}\label{hypcapacity}
\caph K := \lim_{ n \to + \infty} d_{n,h}^{\D} (K) .
\end{equation}

Hyperbolic $n$-th diameter and hyperbolic capacity are invariant under conformal maps. Suppose $f: \D \to U $ is conformal and $K \subset \D$ is compact. Then 
$$ d_{n,h}^U (f(K))= d_{n,h}^{\D}(K) \quad \text{and} \quad \caph_U f(K)= \caph K . $$

\section{Auxiliary Lemmas}\label{lemmas}

The first Lemma stated is a result that derives from elementary calculus. 

\begin{lemma}\label{calculus}
The function $$\sigma(x) := \frac{\log ( -\log \tanh x)}{x}$$ is decreasing in $(0,+ \infty)$ with $\lim_{x \to + \infty} \sigma(x)=-2$. Moreover, there exists a point $x_0$ and a positive constant $c$ such that $$\sigma(x) \leq \frac{\log c}{x}-2, \quad x \geq x_0 . $$
\begin{proof}
With a change of variable, we can see that $\sigma$ has the same kind of monotonicity in $(0,+\infty)$ as the function $$h(x):=\frac{\log ( -\log x)}{\arctanh x}, \quad x \in (0,1).$$
 
The derivative $$h^{\prime }(x) = \frac{   \arctanh x- \frac{x \log x}{1- x^2} \log(-\log x)}{x \log x \arctanh^2 x}.$$
Hence $h^{\prime}$ has the opposite sign of the function
$$ g(x):=\arctanh x - \frac{x \log x}{1- x^2} \log(-\log x), \quad x \in (0,1)$$ considering the fact that $\log x <0$.  
We can easily observe that $g(x)>0$, for $x \in (0, \frac{1}{e}]$. 
We calculate $$ g^{\prime}(x)= - \frac{1}{(1-x^2)^2} \log(- \log x ) \left[(1+x^2) \log x +1-x^2 \right]$$ and we get that $g^{\prime}(x)<0$, for $x \in (\frac{1}{e},1)$. Therefore $g$ is a decreasing function of $x \in (\frac{1}{e}, 1)$ and greater than $$\lim_{x \to 1} g(x)= \lim_{x \to 1} \left[\arctanh x - \frac{x \log x}{1- x^2} \log(-\log x) \right]= \frac{\log 2}{2} >0 . $$
As a result, $h^{\prime}(x) < 0 $, $x \in (0,1)$, which leads to $h$ and $\sigma$ being decreasing functions. 
With elementary calculations and by applying the De L'Hospital rule, the limit 
$$\lim_{x \to + \infty} \frac{\log ( -\log \tanh x)}{x}$$
is equal to $-2$. 
Also from calculus, it can be proved that there exist a positive constant $c$ and a point $x_0 > 0$ such that
$$ \log \frac{e^x+1}{e^x-1} \leq c e^{-x},$$
for every $x \geq x_0$.
The above inequality can easily be interpreted as 
$$- \log \tanh \frac{x}{2} \leq c e^{-x}.$$
Therefore $$
\sigma(x) = h( \tanh x)  = \frac{\log ( -\log \tanh x)}{x}  \leq \frac{\log c -2x}{x} = \frac{\log c}{x}-2,  $$
for every $x \geq x_0$. 
\end{proof}

\end{lemma}

We will also need the hyperbolic distance in the classical Koebe domain. 

\begin{lemma}\label{greenkoebe}
Let $\Delta:= \mathbb{C} \setminus (- \infty, -\frac{1}{4}]$ be the Koebe domain. 
Then for $z,w \in \Delta$, the hyperbolic distance is
\begin{equation}\label{greenkoebeformula}
 d_{\Delta}(z,w) = \arctanh \left| \frac{(4z+1)^{\frac{1}{2}}-(4w+1)^{\frac{1}{2}}}{\overline{(4z+1)}^{\frac{1}{2}}+(4w+1)^{\frac{1}{2}}} \right|.
 \end{equation} 

\begin{proof}
The Koebe function $k: \D \to D $ is written as 
\begin{equation}\label{inversekoebe}
k(z)= \frac{z}{(1-z)^2}= \frac{1}{4} \left( \frac{1+z}{1-z} \right)^2 -\frac{1}{4} \Rightarrow k^{-1}(z) = \frac{\sqrt{4z+1}-1}{\sqrt{4z+1}+1}.
\end{equation}

Due to the conformal invariance of the Green function, we have 
$$d_{\Delta}(z,w) = d_{\D}(k^{-1}(z),k^{-1}(w)) \underset{\eqref{hyperbolicdistance}}{=} \arctanh \left| \frac{k^{-1}(z) -k^{-1}(w) }{1- \overline{k^{-1}(z)}k^{-1}(w)}\right|.$$

Applying \eqref{inversekoebe} and doing elementary calculations, we get 
\eqref{greenkoebeformula}. 

\end{proof}

\end{lemma}

In what follows, we assume that $(\phi_t)$ is a semigroup of holomorphic self-maps of $\D$ and $K$ is a compact non-polar subset of $\D$.

\begin{lemma}\label{limcaph}
The limit $$ \lim_{t \to +\infty}\frac{ \log \caph \phi_t(K)}{t}=0.$$
\begin{proof}

In the case where the semigroup is either hyperbolic or parabolic of positive hyperbolic step, it is proved in \cite{bkkp2} that $\caph \phi_t(K) \xrightarrow{t \to + \infty} \caph_{\bs} h(K) \in \mathbb{R}^{+}$. 
Thus the limit of $\frac{\log \caph \phi_t(K)}{t}$, as $t \to + \infty$, is, in fact, zero. So, we are left to calculate the limit in the case where  the semigroup $(\phi_t)$ is parabolic of zero hyperbolic step. 

We use the associated Koenigs function $h$ of the semigroup $(\phi_t)$ and the associated planar domain $\Omega$. 
Set $\delta := \dist(h(K),\partial \Omega)$. Let $w_1 \in \partial \Omega$ be a point such that $|w_1| = \delta$. 
Due to the fact that $\Omega $ is convex in the horizontal direction, the half-line $\{ w_1 - s : s \geq 0 \} \subset \mathbb{C} \setminus \Omega $. 
Hence $\Omega $ is included into a Koebe-type domain $\Delta_1 := \mathbb{C} \setminus \{ w_1 - s : s \geq 0 \}  $. We can find a linear mapping $f(z) = a z+ b $, with $a>0$ and $ b \in \mathbb{C}$, such that $\Delta_1$ is mapped conformally onto $\Delta := \mathbb{C} \setminus ( - \infty, - \frac{1}{4}]$. Due to the \eqref{subhypden} and the conformal invariance, we obtain
\begin{equation}\label{greater}
\caph \phi_t(K) = \caph_{\Omega} (h(K)+t) \geq \caph_{\Delta_1}(h(K)+t)= \caph_{\Delta} f(h(K)+t). 
\end{equation}
 
Fix $w_1, w_2$ two distinct points in $h(K)$. Set $c:=b+ \frac{1}{4}$ and $c^{\prime} : = \frac{c}{a}$.  The pseudo-hyperbolic distance in $\Delta$ is equal to 
\begin{align*}
\rho_{\Delta}(f(z+t), f(w+t)) &= \rho_{\Delta} ( aw_1+ a t +b, aw_2+at+b) = \left|  \frac{\sqrt{aw_1+at+c}-\sqrt{aw_2+at+c}}{\overline{\sqrt{aw_1+at+c}}+\sqrt{a w_2+at+c}} \right|\\ 
& = |w_1-w_2| \left| \sqrt{w_1+t+ c^{\prime}}+ \sqrt{w_2+t+c^{\prime}} \right|^{-1} \left| \overline{\sqrt{w_1+t+c^{\prime}}}+ \sqrt{w_2+t+c^{\prime}} \right|^{-1}\\ 
&=  |w_1-w_2| \left| |w_1+t +c^{\prime}|+ w_2+ t+c^{\prime} + 2 \sqrt{ w_2+t + c^{\prime}} \RE \{\sqrt{w_1+t+c^{\prime }} \} \right|^{-1} . 
\end{align*}

We can also write 
\begin{align*}
&\left| |w_1+t +c^{\prime}|+ w_2+ t+c^{\prime} + 2 \sqrt{ w_2+t + c^{\prime}}  \RE \{\sqrt{w_1+t+c^{\prime }} \} \right|  =\\
& \quad \quad =   t \left| \left|\frac{w_1+c^{\prime}}{t}+1 \right|+ \frac{w_2+c^{\prime}}{t}+ 1 + 2 \sqrt{ \frac{w_2+ c^{\prime}}{t}+1 } \, \RE \left\{\sqrt{\frac{w_1+c^{\prime }}{t} +1} \right\} \right| \\
& \quad \quad \leq  t \left[ \left|\frac{w_1+c^{\prime}}{t}+1 \right|^{\frac{1}{2}} + \left| \frac{w_2+c^{\prime}}{t}+1 \right|^{\frac{1}{2}} \right]^2  \leq t \left[ 2 \sqrt{\frac{A}{t}+1}\right]^2 \leq 4A + 4t,
\end{align*}
considering the fact that $w_1, w_2$ lie on a compact set and $A:= \max \{|w+ c^{\prime}| : w \in h(K) \}$. 
Hence, 
\begin{equation} \label{ineq1}
\rho_{\Delta}(f(z+t), f(w+t)) \geq \frac{|w_1-w_2|}{4A+4t}. 
\end{equation}

The hyperbolic $n$-th diameter of $h(K)+t$ with respect to $\Delta$ is given by
\begin{align*}
d_{n,h}^{\Delta}(a(h(K)+t)+b)^{\frac{n(n-1)}{2}} & \, = \sup_{w_{\mu}, w_{\nu} \in h(K)+t } \prod_{1 \leq \mu < \nu \leq n} \rho_{\Delta}(a w_{\mu}+b, aw_{\nu}+b) \\
& \, = \sup_{w_{\mu}, w_{\nu} \in h(K) } \prod_{1 \leq \mu < \nu \leq n} \rho_{\Delta}(a w_{\mu}+ at+ b, aw_{\nu}+at +b) \\
& \underset{\eqref{ineq1}}{\geq} (4A+4t)^{- \frac{n(n-1)}{2}} \sup_{w_{\mu}, w_{\nu} \in h(K) } \prod_{1 \leq \mu < \nu \leq n}  |w_{\mu}- w_{\nu}| \\ 
& \underset{\eqref{eudiameter}}{=} (4A+4t)^{- \frac{n(n-1)}{2}} d_n (h(K))^{\frac{n(n-1)}{2}} 
\end{align*}
and hence, applying the logarithm we get
$$\log d_{n,h}^{\Delta}(a(h(K)+t)+b) \geq -\log(A+t)-\log 4 + \log  d_n (h(K)).$$
Taking the limit as $n\to +\infty$, we obtain
$$\log \caph_{\Delta} \left[ a(h(K)+t)+b \right] \geq -\log(A+t) - \log 4+ \log \Capac h(K)$$
and dividing by $t$, 
\begin{equation}\label{ineq2}
\frac{\log \caph_{\Delta} \left[ a(h(K)+t)+b \right]}{t} \geq - \frac{\log(A+t) }{t} - \frac{\log 4}{t} + \frac{\log \Capac h(K)}{t}. 
\end{equation}
Let's note that since $K$ is non-polar, then $h(K)$ is also non-polar (\cite[Corollary 3.6.6]{ransford}) and so, its logarithmic capacity is strictly positive. 
As a result, combining \eqref{greater} and \eqref{ineq2}, we obtain 
\begin{equation}\label{concl1}
\frac{\log \caph \phi_t(K)}{t} \geq - \frac{\log(A+t) }{t} - \frac{\log 4}{t} + \frac{\log \Capac h(K)}{t} \xrightarrow{t \to +\infty} 0. 
\end{equation}

Bearing also in mind that $\caph \phi_t(K)$ is a decreasing function of $t$ and $\caph \phi_t(K) \xrightarrow{t \to + \infty} 0$, as proved in \cite{bkkp2}, it follows that for sufficiently large $t$, $\log \caph \phi_t(K) \leq 0$. 
In conclusion, the limit of $ \frac{\log \caph \phi_t(K)}{t}$, as $t \to + \infty$,  exists and it is equal to $0$.
\end{proof}
\end{lemma}

\begin{lemma}\cite[p.111]{margarnett}\label{111}
Let $E$ be a compact non-polar subset of $\D$ and let $\nu $ be its Green equilibrium measure. 
For any $z \in \D \setminus E$,
$$ \omega (z, E, \D) =\caph E \cdot G_{\nu} (z)   .$$
\end{lemma}

\begin{lemma}\label{inequality}
For every $ z \in K$, there exists a constant $\widetilde{C}$ such that 
 $$ \omega(z, \phi_t(K), \D) \leq \widetilde{C} e^{- \pi \lambda_{\D}(K, \phi_t(K)) } .$$
\begin{proof}
Suppose $z \in K$. Consider the multiply connected domain $D:=\D \setminus \phi_t(K)$. From \eqref{generh7}, we have that $$\omega (z, \phi_t(K), \D) \leq \frac{8}{\pi} \exp \left\{ - \pi \sup_{ \sigma} \lambda_{D} (\sigma, \phi_t(K)) \right\}, $$
where $\sigma$ is any Jordan arc in $\D$ that joins $z$ with $\partial \D$. 
As $t$ increases, $\phi_t(K)$ approaches the Denjoy-Wolff point $1$ and as $t \to + \infty$, it shrinks to a point. There exists $t_1 \in \mathbb{R}^{+}$ such that for all $t \geq t_1$, $K \cap \phi_t(K) \neq \emptyset$. 
Pick $z_0 \in \partial K$ such that it minimizes the distance from $\phi_t(K)$.
Choose two arcs $\sigma_1$ and $\sigma_2$ in $\D \setminus (K \cup \phi_t(K))$, pairwise disjoint, which join $z_0$ with $\partial \D$. We further suppose that $\gamma := \sigma_1 \cup \sigma_2$ is a crosscut and $K$ and $\phi_t(K)$ lie on different components of $\D \setminus \gamma$. 

For $i=1,2$, $$\lambda_{D}(\sigma_i,\phi_t(K))\leq \sup_{\sigma} \lambda_D (\sigma, \phi_t(K))$$ and since the extremal distance decreases as any of the sets increase (see \cite[Corollary p.54]{ahlforsconformal}), we obtain $$\lambda_{D} (\gamma, \phi_t(K)) = \lambda_{D}(\sigma_1 \cup \sigma_2,\phi_t(K))\leq \sup_{\sigma} \lambda_D (\sigma, \phi_t(K)). $$ 

\begin{figure}[H]
\begin{center}
\includegraphics[scale=0.29]{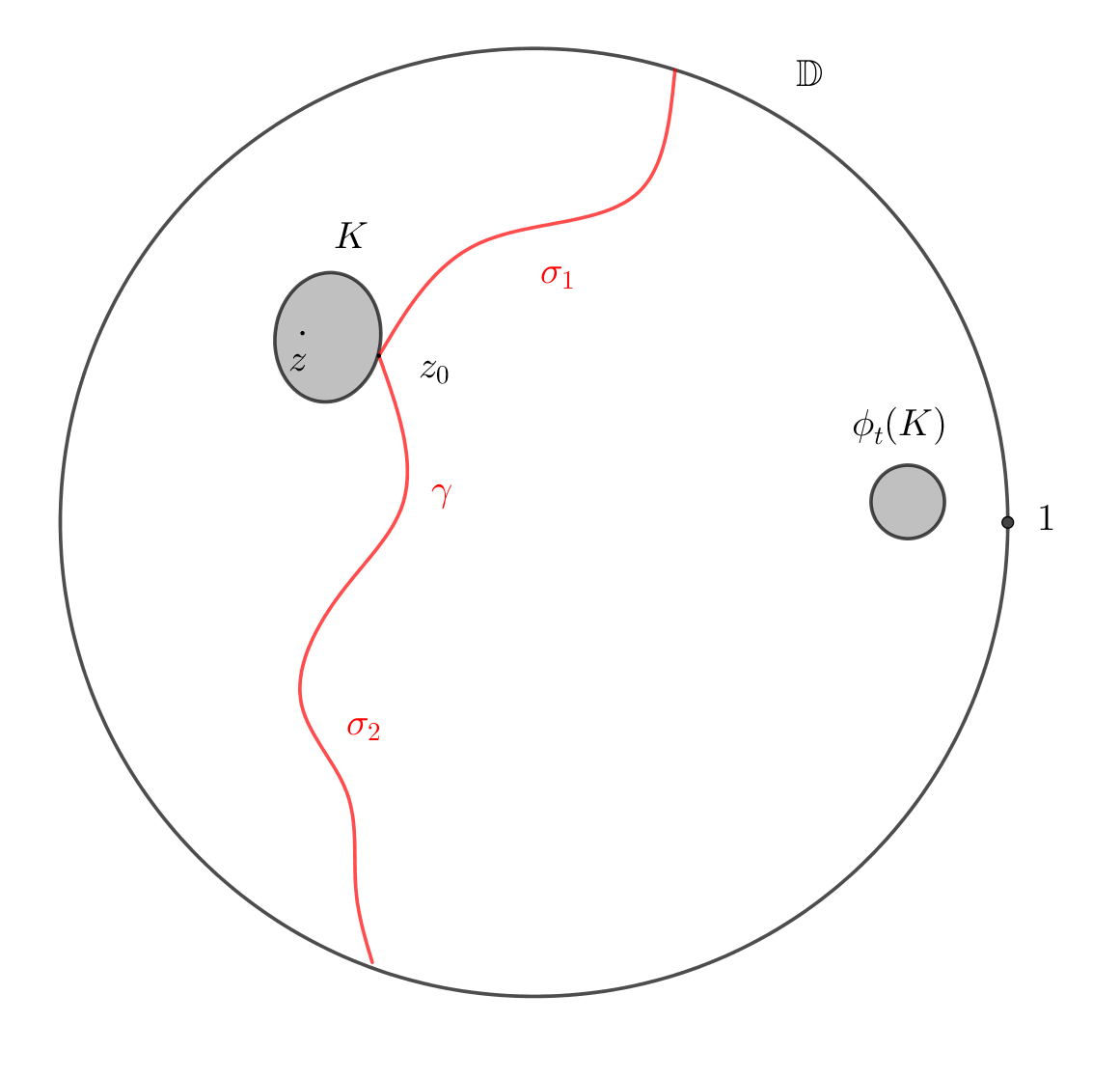}
\caption{Arcs $\sigma$}
\end{center}
\end{figure} 
 
Let $ \{ \Gamma^{\star} \}$ denote the family of rectifiable curves that separate $K$ from $\phi_t(K)$ and $ \{ \Gamma_{\gamma}^{\star} \}$ the family of rectifiable curves that separate $\gamma$ from $\phi_t(K)$. Let $\{ \Gamma_1^{\star} \} $ denote the family of rectifiable curves that separate $K$ from $\phi_t(K)$ but not $\phi_t(K)$ from $\gamma$.  
Obviously, $\{ \Gamma^{\star} \}$ is the disjoint union $\{ \Gamma_{\gamma}^{\star} \} \cup \{ \Gamma_1^{\star} \} $. 
From Lemma \ref{non-disjointparallel} and the definition of the conjugate extremal length, 
$$\lambda_{D} ( K,\phi_t(K)) \leq \lambda_{D} (\gamma, \phi_t(K)) +\frac{1}{\lambda_D (\{ \Gamma_1^{\star} \}) } \leq \sup_{\sigma} \lambda_D (\sigma, \phi_t(K)) +\frac{1}{\lambda_D (\{ \Gamma_1^{\star} \}) }. $$

If the family $ \{ \Gamma_1^{\star} \}$ is empty, then $\lambda_D (\{ \Gamma_1^{\star} \})  = + \infty$ (\cite[p.70]{ohtsuka}) and the above inequality holds as equality. 
In the case where $\lambda_D (\{ \Gamma_1^{\star}\})=0$, according to the definition, we get $$ \sup_{\rho} \frac{ \inf_{\gamma \in \{ \Gamma_1^{\star} \}} \Le( \gamma, \rho)^2 }{\Ar(\D, \rho)} =0 \Rightarrow \inf_{\gamma \in \{ \Gamma_1^{\star} \}} \Le( \gamma, \rho) =0 , $$
for any measurable function $\rho$. 
The same holds for the euclidean metric on the unit disk, as $e(z)=1$ is a non-negative Borel measurable function on $\D$ with finite area $\Ar ( e, \D)= \pi$. Hence 
$$\Le (\{\Gamma_1^{\star} \}, e)=0  \Leftrightarrow \inf_{\gamma \in \{\Gamma_1^{\star} \}} \int_{ \gamma} |dz| =0. $$ 
Therefore for every $\epsilon > 0 $, by the definition of the infimum, there exists at least one curve $\gamma \in \{ \Gamma_1^{\star} \} $ such that 
$\int_{\gamma} |dz| <\epsilon. $
However such a curve cannot separate the compact sets $K$ and $\phi_t(K)$. Hence $ \lambda_D(\{ \Gamma_1^{\star} \})>0$. All the curves that separate $K$ and $\phi_t(K)$ have strictly positive euclidean length and so $ \lambda_D(\{ \Gamma_1^{\star} \})$ has a positive lower bound. 
As a result, there exists a constant $C_1>0$ such that $\lambda_D (\{ \Gamma_1^{\star}\})  \geq C_1$ and  
$$ \lambda_{\D} ( K,\phi_t(K))=\lambda_D( K,\phi_t(K)) \leq  \sup_{\sigma} \lambda_D (\sigma, \phi_t(K))+ \frac{1}{C_1} ,$$
where the first equality follows from the definition of extremal distance; see \eqref{extconj}. 
We obtain 
\begin{equation}
\omega (z_0, \phi_t(K), \D) \leq \frac{8 }{\pi} \exp \left\{- \pi \lambda_{\D}(K,\phi_t(K))+\frac{\pi}{C_1} \right\}.
\end{equation}
According to Harnack's inequality (see e.g. \cite[Definition 1.3.4]{ransford} and \cite[p.40]{margarnett}), 
$$\tau^{-1}(z,z_0) \cdot \omega (z, \phi_t(K), \D) \leq \omega (z_0, \phi_t(K), \D) \leq \frac{8}{\pi} \exp \left\{- \pi \lambda_{\D}(K,\phi_t(K))+\frac{\pi}{C_1} \right\},$$
where $\tau$ denotes the Harnack distance in $\D \setminus \phi_t(K)$.
Set $c:= \sup_{w \in K} \tau(w,z_0) > 0$, which exists as maximum, since $K$ is compact and $\tau$ is continuous. 
Hence 
\begin{equation}\label{upperbound}
\omega (z, \phi_t(K), \D) \leq \frac{8 c}{ \pi} \exp \left\{- \pi \lambda_{\D}(K,\phi_t(K))+\frac{\pi}{C_1} \right\}  = :\widetilde{C} e^{- \pi \lambda_{\D}(K,\phi_t(K))} .
\end{equation}
\end{proof}
\end{lemma}

\section{Proof of Theorem \ref{main}}\label{proof}

In what follows, we will prove Theorem \ref{main} and the convergence to the spectral value for every potential theoretic quantity.

\subsection{Hyperbolic Metric}

Fix $z \in K$. The limit of $\frac{d_{\D}(z, \phi_t(z))}{t}$, as $t \to + \infty$, exists and it is equal to $\frac{\lambda}{2}$, as seen in \eqref{divrate}. 
Since $d_{\D}$ is a metric
\begin{align*}
\sup_{w \in K} \left| \frac{d_{\D}(z, \phi_t(w))}{t} - \frac{\lambda}{2} \right| & \leq \sup \bigg\{ \sup_{w \in K} \left| \frac{d_{\D}(z, \phi_t (z))}{t} - \frac{d_{\D}(\phi_t(z), \phi_t(w))}{t} -\frac{\lambda}{2} \right|,  \\
& \quad \quad \quad   \sup_{w \in K}  \left|  \frac{d_{\D}(z, \phi_t (z))}{t}+  \frac{d_{\D}(\phi_t(z), \phi_t(w))}{t} -\frac{\lambda}{2} \right| \bigg\}  \\
& \leq \left|\frac{d_{\D}(z, \phi_t (z))}{t} -\frac{\lambda}{2} \right| + \sup_{w \in K} \left| \frac{d_{\D}(\phi_t(z), \phi_t(w))}{t} \right| \\
& < \left|\frac{d_{\D}(z, \phi_t (z))}{t} -\frac{\lambda}{2} \right|+ \sup_{w \in K} \left| \frac{d_{\D}(z,w)}{t} \right|,
\end{align*}
where the last inequality follows from Schwarz-Pick lemma. 
Since $z,w$ lie on a compact set, there exists a constant $c>0$ such that their hyperbolic distance $d_{\D}(z,w) \leq c$. Therefore, 
\begin{equation}\label{limhypmet}
\sup_{w \in K} \left| \frac{d_{\D}(z, \phi_t(w))}{t} - \frac{\lambda}{2} \right| < \left|\frac{d_{\D}(z, \phi_t (z))}{t} -\frac{\lambda}{2} \right|+ \frac{c}{t} \xrightarrow{ t \to + \infty} 0 
\end{equation}
and so, $\frac{d_{\D}(z, \phi_t(w))}{t}$ converges uniformly in $w \in K$ to $\frac{\lambda}{2}$, as $t \to + \infty$. \qed

\subsection{Green Function}

Fix $z \in K$. For $w \in K$, set $$u_t(w) : = \frac{\log g_{\D} (z, \phi_t(w))}{ d_{\D} (z, \phi_t(w))} \underset{\eqref{greendist}}{=} \frac{\log [-\log \tanh d_{\D}(z, \phi_t(w))]}{d_{\D}(z, \phi_t(w))}.  $$
According to Lemma \ref{calculus}, the limit 
$$ \lim_{t \to +\infty} u_t(w)  = -2$$
and hence $u_t$ converges pointwise to $-2$ on $K$.
We examine whether the convergence is also uniform. 
Due to Lemma \ref{calculus}, we obtain
$$ 0 \leq u_t(w)+2 \leq \frac{\log c}{d_{\D}(z, \phi_t(w))}$$
and so,  
$$   \sup_{w \in K} |u_t(w)+2 |  \leq \sup_{w \in K}  \left| \frac{\log c}{d_{\D}(z, \phi_t(w))} \right|  \leq  \frac{|\log c|}{\inf_{w \in h(K)} d_{\D}(z, \phi_t(w))} \xrightarrow{t \to + \infty} 0 ,$$
since $\phi_t(w) \xrightarrow{t \to + \infty} 1$, for every $w \in K$. As a result, $u_t$ converges uniformly in $w \in K$ to $-2$, as $t \to + \infty$. 

For $t$ sufficiently large, we can verify that $d_{\D}(z, \phi_t(w)) > const. >0$, for any $w \in K$, both $\frac{d_{\D}(z, \phi_t(w))}{t}$ and $u_t(w)$ are bounded on $K$, and since they converge uniformly on $K$, their product also converges uniformly on $K$. 
Therefore,  
\begin{equation}\label{limgreenfunction}
 \frac{\log g_{\D}(z, \phi_t(w))}{t} \xrightarrow{t \to + \infty} - \lambda
\end{equation}
uniformly in $w \in K$. \qed

\subsection{ Harmonic Measure}

Fix $t \in \mathbb{R}^{+}$ sufficiently large. 
Since $\phi_t (K)$ is compact and non-polar, it has a unique Green equilibrium measure $\mu_t$. Fix $z \in K$. 
According to Lemma \ref{111},
\begin{equation}\label{eq111}
 \frac{\omega (z, \phi_t(K), \D)}{\caph \phi_t(K)} =  G_{ \mu_t} (z) = \int_{\phi_t(K)} g_{\D} (z, w) d \mu_t (w),
 \end{equation} 
where $g_{\D}$ is the Green function in the unit disk. 
Estimating the above integral, we obtain the following inequality
\begin{equation}\label{estimates}
 \mu_t(\phi_t(K))\inf_{w \in K} g_{\D}(z,\phi_t(w)) \leq \frac{\omega (z, \phi_t(K), \D)}{\caph \phi_t(K) } \leq \mu_t(\phi_t(K)) \sup_{w \in K} g_{\D}(z,\phi_t(w)).
\end{equation}
However, $\mu_t$ is a probability measure and so, $\mu_t(\phi_t(K))=1$. 
Applying the logarithm, which is both continuous and increasing, we get 
\begin{equation}\label{min,max}
\frac{\log \caph \phi_t(K)}{t}  +\inf_{w \in K}   \frac{\log g_{\D}(z,\phi_t(w))}{t} \leq \frac{\log \omega (z, \phi_t(K), \D)}{t}  \leq  \frac{\log \caph \phi_t(K)}{t} + \sup_{w \in K} \frac{\log g_{\D}(z,\phi_t(w))}{t}.
\end{equation}
From \eqref{limgreenfunction} and Lemma \ref{limcaph}, we conclude that 
\begin{equation}\label{limha}
\lim_{t \to + \infty} \frac{\log \omega (z, \phi_t(K), \D)}{t} = -\lambda.
\end{equation}
\qed

\subsection{Extremal Length }

From Lemma \ref{inequality}, we obtain 
\begin{equation}\label{upperbound1}
\pi \lambda_{\D}(K,\phi_t(K)) \leq \log \widetilde{C} - \log \omega(z, \phi_t(K),\D) \Rightarrow \pi \frac{\lambda_{\D}(K,\phi_t(K))}{t} \leq \frac{\log \widetilde{C}}{t} - \frac{\log \omega(z, \phi_t(K),\D)}{t}.
\end{equation}
Hence
\begin{equation}\label{ineqfinalext1}
 \limsup_{t \to + \infty} \frac{\lambda_{\D}(K,\phi_t(K))}{t} \leq \frac{\lambda}{\pi} . 
\end{equation}

In order to find a lower bound, we use a technique with crosscuts from \cite[\S\ 2.19]{ohtsuka}. Let $h$ be the associated Koenigs function of the semigroup. Due to its linearization property, $h(\phi_t(K))= h(K)+t$. 
Suppose $y:=\frac{1}{2} \left(\min \{\IM w :w \in h(K)\} + \max \{\IM w : w \in h(K)\}  \right)$ and take the horizontal line $L:= \{ z :  \IM z=y \}$. If $L \cap \partial \Omega \neq \emptyset$, fix $w \in \partial \Omega \cap L$ such that $\RE w = \max\{ \RE \widetilde{w} : \widetilde{w} \in \partial \Omega \cap L \} $.  
In the case where $L \cap \partial \Omega = \emptyset$, fix $w \in \partial \Omega$ such that $\dist(w, h(K))= \dist (\partial \Omega, h(K))$; see Figures \ref{capnonempty} and \ref{capempty}. 
Set $$z_t:= i y+ \frac{t+ \min_{\widetilde{w} \in h(K)} \RE \widetilde{w} + \max_{\widetilde{w} \in h(K)} \RE \widetilde{w}}{2} \in L$$ and $R_t:=|w -z_t|$. 

\begin{minipage}{\linewidth}
   \centering
     \begin{minipage}{0.45\linewidth}
     \begin{figure}[H]
     \includegraphics[scale=0.4]{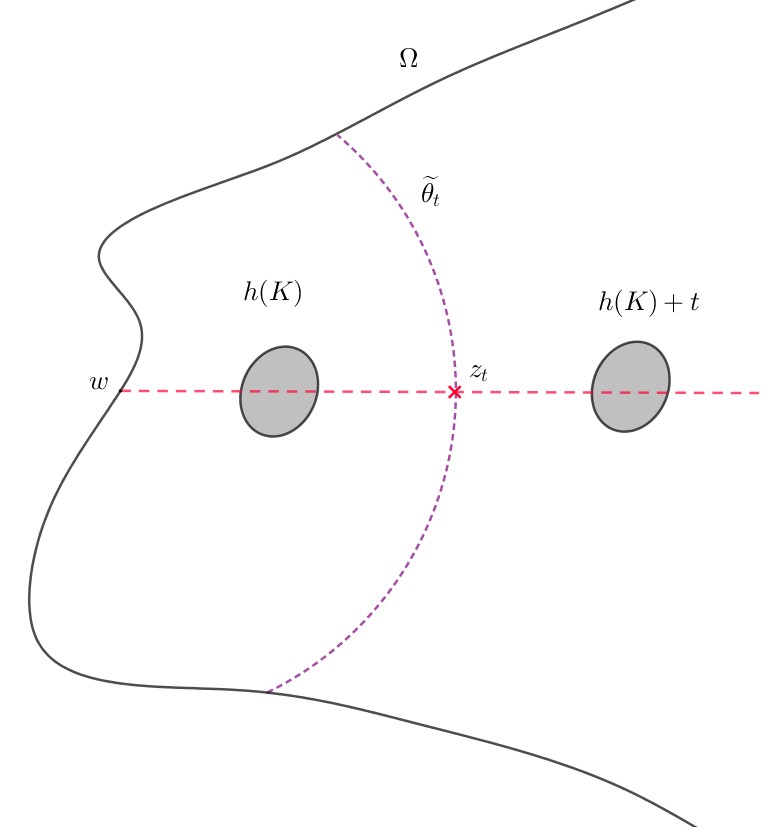}
     \caption{$L \cap \partial \Omega \neq \emptyset$}
     \label{capnonempty}
     \end{figure}
      \end{minipage}
      \hspace{0.05\linewidth}
      \begin{minipage}{0.45\linewidth}
     \begin{figure}[H]
     \includegraphics[scale=0.4]{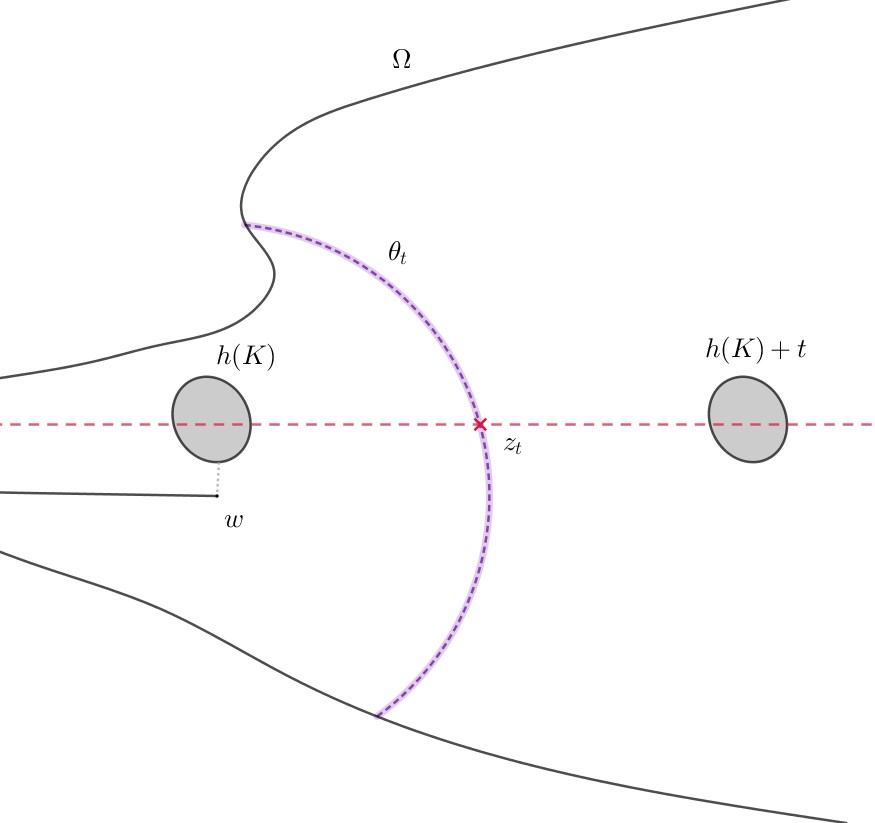}
     \caption{$L \cap \partial \Omega = \emptyset$}
     \label{capempty}
    \end{figure}
   \end{minipage}
 \end{minipage}
\vspace*{0.2cm}

Fix $t_0 := \inf \{ t \in \mathbb{R}^{+}: h(K) \subset D(w, R_t)\}$. The existence of such $t_0$ is guaranteed by the fact that $h(K)$ is compact and $R_t$ is increasing. 
Moreover, for some $t_1>0$ onward, $D(w,R_t) \cap (h(K)+t )= \emptyset$. Take $\widetilde{t_0} = \max\{ t_0,t_1\}$. 
Define $\widetilde{\theta_0}$ to be the crosscut in $\Omega $ that is the connected component of $\partial D(w, R_{\widetilde{t_0}}) \cap \Omega$, which passes through $z_{\widetilde{t_0}}$. 
Then $\widetilde{\theta_0}$ separates $h(K)$ from $h(K)+\widetilde{t_0}$. 

\begin{figure}
\begin{center}
\includegraphics[scale=0.45]{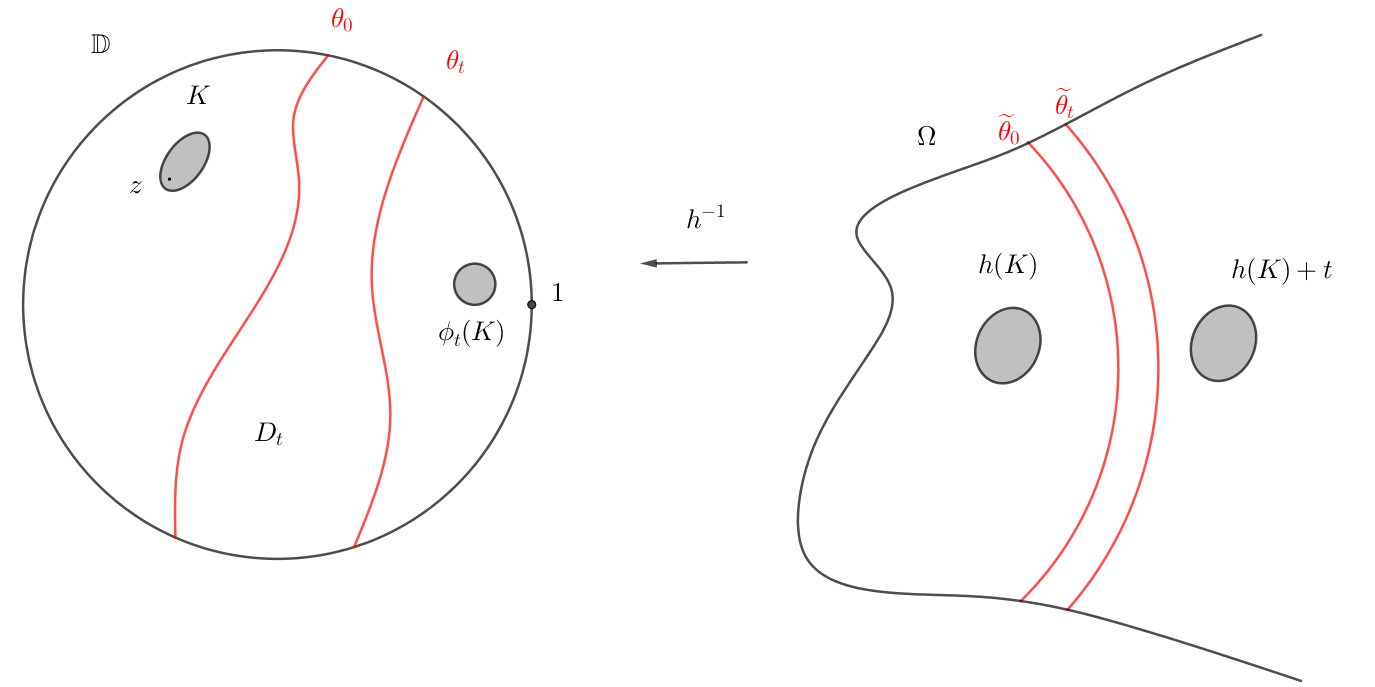}
\caption{Families of crosscuts $\{\theta_t \}$ and $\{ \widetilde{\theta}_t \}$}
\label{cross}
\end{center}
\end{figure}

In the same way, we define a family of crosscuts $\{ \widetilde{ \theta_t} \}_{t \geq \widetilde{t_0}} $ in $\Omega$, each of which is the connected component of $\partial D(w, R_t) \cap \Omega$ that passes through $z_t$ and separates $h(K)$ from $h(K)+t$. As $t \to + \infty$, $z_t \to \infty$ and $R_t \to + \infty$. Hence the family of the disks $\{D(w, R_t)\}_t$ is increasing and covers the whole complex plane. 
Let's also note that any two crosscuts of the family have distinct end-points and are pairwise disjoint. 
Taking their pre-images with respect to the Koenigs function, we obtain a family of crosscuts $ \{\theta_t \}_t$ in $\D$, such that for each $t$, $\theta_t:= h^{-1} ( \widetilde{\theta_t} )$; see Figure \ref{cross}.
From \cite[Prop. 2.14]{pommerenke4}, the crosscuts $\theta_t$ have distinct endpoints on the unit circle and each separates $K$ from $\phi_t(K)$. 

Due to connectedness, the domains $D_t:=h^{-1}(D(w,R_t) \cap \Omega)$ form an increasing family of domains that converges to $\D$, as $t \to + \infty$. 
Moreover, as $t \to + \infty$, the euclidean diameter of $\theta_t$ tends to $0 $, so, $\{ \theta_t\}$ forms a null-chain that separates $z$ from  the the Denjoy-Wolff point of the semigroup. 
For the detailed theory on crosscuts and null-chains, the reader may refer to \cite[\S\ 2.4]{pommerenke4}.

Set $\theta_0 := h^{-1} ( \widetilde{\theta}_0 )$. From \cite[Theorem 2.82]{ohtsuka}, it is true that 
\begin{equation}\label{2.82ohtsuka}
\omega (z, \theta_t, D_t )= \exp \left\{ - \pi \lambda_{\D} (\theta_0, \theta_t) + \mathcal{O}(1) \right\}, 
\end{equation}
as $t \to + \infty$. According to \cite[Theorem III.5.1 and (5.3)]{margarnett} and since $\{D_t\}$ is an increasing family of domains, we get 
\begin{equation}\label{lim1}
\lim_{t \to + \infty} \omega(z, \theta_t, D_t) = \omega (z, \{1\}, \D) = 0, \quad z \in K,
\end{equation}
as $\{1\}$ is a polar subset of $\partial \D$. 
Moreover, for $z \in K$ 
$$\lim_{t \to + \infty} \omega(z, \phi_t(K), \D)  = \lim_{t \to +\infty} \omega (h(z), h(K)+t, \Omega) = \omega( h(z), \{ \infty \} , \mathbb{C})=0, $$
as $\{ \infty \}$ is polar. 

From Proposition \ref{markov}, we observe that $\omega (h(z), \widetilde{\theta_{t}}, D(w,R_t) \cap \Omega )$ is a decreasing function of $t$. 
Indeed, for $ \widetilde{t_0}  < t_1 < t_2$, we obtain 
\begin{align*}
\omega (h(z), \widetilde{\theta_{t_2}}, D(w,R_{t_2}) \cap \Omega) & =  \int_{\widetilde{ \theta_{t_1}}} \omega (h(z), ds,  D(w,R_{t_1}) \cap \Omega) \cdot \omega (s,\widetilde{ \theta_{t_2}} ,D(w,R_{t_2}) \cap \Omega) \\
& \leq \omega (h(z), \widetilde{ \theta_{t_1}},  D(w,R_{t_1}) \cap \Omega),
\end{align*}
where the inequality follows due to the fact that the harmonic measure is a probability measure. 
Furthermore, using the probabilistic interpretation of harmonic measure (see Subsection \ref{harmonicmeasure}), we can observe that $ \omega(h(z), h(K)+t, \Omega)$ is also a decreasing function of $t$, since $h(K)+t$ moves towards infinity, as $t$ increases. Due to conformal invariance, it follows that both $\omega(z, \theta_t, D_t)$ and $\omega(z, \phi_t(K), \D)$ are decreasing functions of $t$ and hence, they satisfy Dini's criterion for uniform convergence on the compact set $K$. So, for every  $ \epsilon >0 $, 
\begin{equation}\label{smallo}
\sup_{z \in K} | \omega(z, \phi_t(K), \D)  - \omega(z, \theta_t, D_t) |<  \epsilon .
\end{equation}

Let's denote by $\{ \Gamma_t^{\star}\}$ the family of curves that separate $\theta_0$ from $\theta_t$ and by $\{ \Gamma_K^{\star}\} $ the family of curves that separate $K$ from $ \phi_t(K)$. Obviously $\{ \Gamma_t^{\star} \} \subset \{ \Gamma_K^{\star} \}$. From Lemma \ref{extension}, $\lambda (\{ \Gamma_t^{\star} \}) \geq \lambda(\{ \Gamma_K^{\star} \})$, which leads to 
$\lambda_{\D}^{\star} ( \theta_0, \theta_t) \geq \lambda_{\D}^{\star} (K, \phi_t(K))$. Therefore, as expected, 
\begin{equation}\label{extdist}
\lambda_{\D}( \theta_0, \theta_t) \leq \lambda_{\D} (K, \phi_t(K)).
\end{equation}

Applying \eqref{smallo} to \eqref{2.82ohtsuka}, we get
\begin{equation} \label{lowerbound}
\omega (z, \phi_t(K), \D) > \omega(z, \theta_t, D_t) -  \epsilon =\exp \{ - \pi \lambda_{\D} ( \theta_0, \theta_t) + \mathcal{O}(1) \} -  \epsilon, 
\end{equation}
for sufficiently large $t$. Also, for large enough $t$, $\mathcal{O}(1)$ is bounded from below by a constant $ \log C$, where $C>0$. Therefore, from \eqref{extdist} and \eqref{lowerbound}, 
$$
\omega (z, \phi_t(K), \D) \geq C  e^{ - \pi \lambda_{\D} (K, \phi_t(K)) }  -  \epsilon \Rightarrow \frac{\log \omega (z, \phi_t(K), \D)}{t} > \frac{ \log \left( C e^{-\pi \lambda_{\D} (K, \phi_t(K))} -  \epsilon \right)}{t}.$$
Using \eqref{limha}, we obtain for every $\epsilon >0$ and every sufficiently large $t$,
$$-\lambda+  \epsilon > \frac{ \log \left( C e^{-\pi \lambda_{\D} (K, \phi_t(K))} - \epsilon \right)}{t} \Rightarrow \log \frac{\epsilon + e^{ \epsilon t} e^{-\lambda t}}{C} > - \pi \lambda_{\D}(K, \phi_t(K)). $$
Dividing again by $t$, it follows 
\begin{align}\label{ineqext1}
\frac{\lambda_{\D}(K, \phi_t(K))}{t} &> -\frac{1}{\pi t} \log \left(  \epsilon + e^{ \epsilon t} e^{-\lambda t } \right)  + \frac{\log C}{\pi t} \nonumber \\ 
& =-\frac{\epsilon}{\pi} - \frac{1}{\pi t} \log ( \epsilon e^{- \epsilon t } +  e^{-\lambda t })  + \frac{\log C}{\pi t}.  \nonumber 
\end{align}
Note that for sufficiently large $t$, $\frac{1}{t}< \epsilon$ and so, $ \epsilon e^{-\epsilon t }< \epsilon e$. Choose $\epsilon < \frac{1}{e}$. Hence for every $\eta>0$,
$$\frac{\lambda_{\D}(K, \phi_t(K))}{t} >  -\frac{\epsilon}{\pi} - \frac{1}{\pi t} \log ( 1 +  e^{-\lambda t}  )  + \frac{\log C}{\pi t}> \frac{\lambda }{\pi}  - \eta, $$
as $ t^{-1} \log ( 1 +  e^{-\lambda t} ) \xrightarrow{t \to + \infty} -\lambda$. Therefore, we are led to 
\begin{equation}\label{ineqfinalext2}
\liminf_{t \to + \infty} \frac{\lambda_{\D}(K, \phi_t(K))}{t}  \geq \frac{\lambda}{\pi}. 
\end{equation}

As a result, combining \eqref{ineqfinalext1} and \eqref{ineqfinalext2}, we get
$\pi \frac{\lambda_{\D}(K,\phi_t(K))}{t} \to \lambda.$ \qed

\subsection{Condenser Capacity}
 From \eqref{extdistcap}, we obtain $ \frac{1}{t \Capac(\D, K, \phi_t(K))} \xrightarrow{t \to + \infty} \frac{\lambda }{\pi } .$  
If the spectral value is positive, then we can write that 
 \begin{equation}\label{limitcondencap}
 \lim_{t \to + \infty} t \Capac(\D, K, \phi_t(K)) = \frac{\pi}{ \lambda }.
 \end{equation}
In the case the spectral value is zero, then $ \frac{1}{t \Capac(\D, K, \phi_t(K))} \xrightarrow{t \to + \infty} 0 $ and since the condenser capacity is positive, we get $t \Capac(\D, K, \phi_t(K)) \xrightarrow{t \to + \infty} + \infty.$ Hence, the limit \eqref{limitcondencap} holds in the generalized form.   \qed

\section{Proof of Proposition \ref{asymptoticintegral}} \label{proofprop}

Combining Theorem \ref{main} with \eqref{connection}, we obtain that
$$ \frac{\pi}{2} \lim_{u \to + \infty} \frac{\lambda_{\D}(K, \phi_u(K))}{u}=  \lim_{u \to +\infty} \frac{s_u(\phi_t,z)}{u},$$ 
which leads to $  \lambda_{\D}(K, \phi_u(K))= \frac{2}{\pi} s_u(\phi_t,z) + \smallO(u)$, as $u \to +\infty$, 
for any point $z \in \D$.
So, asymptotically the extremal distance between compact sets plays the role of the hyperbolic step. \qed

\section*{Acknowledgments}

The author wishes to thank the Department of Mathematics, University of Wuerzburg, where part of this research was carried out.

This research is co-financed by Greece and the European Union (European Social Fund- ESF) through the Operational Programme ``Human Resources Development, Education and Lifelong Learning'' in the context of the project ``Reinforcement of Post-doctoral researchers - 2nd Cycle'' (MIS-5033021), implemented by the State Scholarships Foundation (IKY). 

\medskip

\bibliographystyle{plain}

\end{document}